\UseRawInputEncoding
\documentclass{article}

\usepackage{stmaryrd}
\usepackage{lmodern}
\usepackage{diagbox}
\usepackage{mathtools}
\usepackage{amssymb,amsmath}
\usepackage{graphicx}

\usepackage{mathrsfs}

\usepackage{color}
\usepackage{amsfonts}
\usepackage{epsfig}
\usepackage{graphics}

\usepackage{multicol}
\usepackage{lscape}
\usepackage{latexsym}
\def\proof{\noindent {\bf Proof }}
\def\Fq{{\mathbb F}_q}

\def\Zq-1{{\mathbb Z}_{q-1}}

\def\qed{~~\vrule height8pt width4pt depth0pt}

\def\Mcal{{\cal M}}
\def\Ical{{\cal I}}

\def\Pcal{{\cal P}}
\def\F2{{\mathbb F}_2}

\def\eps{{\varepsilon}}

\def\Mcal{{\cal M}}
\def\Ecal{{\cal E}}

\newtheorem{prop}{Proposition}

\newtheorem{lemma}{Lemma}
\newtheorem{thm}{Theorem}
\newtheorem{cor}{Corollary}
\newtheorem{example}{Example}

\newcommand{\fq}{{\mathbb F}_{q}}

\newcommand{\fqx}{{\mathbb F}_{q}^{*}}

\newcommand{\Scal}{{\cal S}}

\newcommand{\om}{{\omega}}

\begin{document}

\title {Improved error bounds for the number of irreducible polynomials and self-reciprocal irreducible monic polynomials with prescribed  coefficients over a finite field }
 \author{
Zhicheng Gao\\
School of Mathematics and Statistics\\
Carleton University\\
Ottawa, Ontario\\
Canada K1S5B6\\
Email:~zgao@math.carleton.ca }
\maketitle

\begin{abstract}
A polynomial is called self-reciprocal (or palindromic) if the sequence of its coefficients is palindromic. In this paper we obtain improved error bounds for the number of irreducible polynomials and self-reciprocal irreducible monic polynomials with prescribed  coefficients over a finite field.
The improved bounds imply
self-reciprocal irreducible monic polynomials with degree $2d$ and  prescribed $\ell$ leading coefficients always exist provided that $\ell$ is slightly less than $d/2$.
\end{abstract}

\section{ Introduction}

The main objective of this paper is to derive improved error bounds on the number of irreducible polynomials and  self-reciprocal irreducible monic polynomials  with prescribed coefficients. Practical error bounds played important role in proving the existence of irreducible polynomials and self-reciprocal irreducible monic polynomials with prescribed coefficients. For example, the famous Hansen-Mullen conjecture on the existence of irreducible polynomials with one prescribed coefficient was proved (for large values of $q,d$) by Wan \cite{Wan97} using Weil's bound on character sums. Panario and Tzanakis \cite{PanTza12} used Wan's approach to study the extended Hansen-Mullen conjecture by considering several prescribed coefficients. Ha \cite{Ha16} obtained bounds for several prescribed coefficients using a different approach. Garefalakis and  Kapetanakis \cite{GarKap12} used Wan's approach to prove the existence of self-reciprocal irreducible monic polynomials with one prescribed coefficient. Pollack \cite{Pol13} discussed asymptotic number of irreducible polynomials with several prescribed coefficients.

Throughout the paper, we shall use the  notations from \cite{Gao21}, which are summarized below.
\begin{itemize}
\item $\fq$ denotes the finite field with $q$ elements, where $q=p^r$ for some prime $p$ and positive integer $r$.
\item $\Mcal_q$ denotes the set of monic polynomials over $\Fq$.
\item $\left[x^j\right]f(x)$ denotes the coefficient of $x^j$ in the polynomial $f(x)$.
\item For a polynomial $f$, $\deg(f)$ denotes the degree of $f$, and \\
$\displaystyle f^*(x)=x^{\deg(f)}f(1/x)$ is the {\em reciprocal} of $f$.
\item $\Pcal_q$ denotes the set of polynomials in $\Mcal_q$ with $f^*=f$. Polynomials in $\Pcal_q$ are called {\em self-reciprocal} or {\em palindromic}.
\item $\Ical_q\subseteq \Mcal_q$  denotes the set of irreducible monic polynomials.
\item $\Scal_q=\Ical_q\cap \Pcal_q$ denotes the set of  self-reciprocal  irreducible monic polynomials  over $\Fq$.
\item $\Scal_q(d)=\{f: f\in \Scal_q, \deg(f)=d\}$.
\item $\Ical_q(d)=\{f: f\in \Ical_q, \deg(f)=d\}$.
\end{itemize}
Given non-negative integers $\ell,t$, and vectors $\vec{a}=(a_1,\ldots,a_{\ell})$ and $\vec{b}=(b_0,b_1,\ldots,b_{t-1})$, we define
\begin{align*}
\Ical_q(d;\vec{a},\vec{b})&=\{f:f\in \Ical_q(d), [x^{d-j}]f(x)=a_{j}, 1\le j\le \ell, [x^j]f(x)=b_j, 0\le j\le t-1\},\\
\Scal_q(d;\vec{a})&=\{f:f\in \Scal_q(d), [x^{d-j}]f(x)=a_{j}, 1\le j\le \ell\}.
\end{align*}
Thus the vector $\vec{a}$ gives the $\ell$ {\em leading coefficients} of $f$, and $a_j$ is also called the $j$th {\em trace} of $f$. The vector $\vec{b}$ gives the $t$ {\em ending coefficients} of $f$ and $b_0$ is the {\em norm} of $f$. We shall use $I_q(d;\vec{a},\vec{b}):=\left|\Ical_q(d;\vec{a},\vec{b})\right|$ to denote the cardinality of $\Ical_q(d;\vec{a},\vec{b})$, and $S_q(d;\vec{a})$ to denote the cardinality of $\Scal_q(2d;\vec{a})$.

It is clear $I_q(d;\vec{a},\vec{b})=0$ when $d>1$ and $b_0=0$. Hence we assume $b_0\ne 0$ when $t\ge 1$.
It is also clear $S_q(d/2;\vec{a})=0$ when $d>1$ is odd.

The rest of the paper is organized as follows. In section~2 we state our main results about the bounds for $I_q(d;\eps)$ and $S_q(d;\eps)$. Proofs are given in section~3.
In section~4, we give some examples to demonstrate the improvement of our bounds over those in \cite{Coh05,Hsu96}. Section~5 concludes the paper.

\section{Main results}
We need a few more notations before stating our main results. Following Hayes \cite{Hay65},
we say that two polynomials $f,g\in \Mcal$ are {\em equivalent} with respect to $\ell,t$ if
\begin{align*}
\left[x^{\deg(f)-j}\right]f(x)&=\left[x^{\deg(g)-j}\right]g(x), 1\le j\le \ell, \\
\left[x^j\right]f(x)&=\left[x^j\right]g(x), 0\le j\le t-1.
\end{align*}
Let $\langle f\rangle$ denote the equivalence class represented by $f$.  It is known  \cite{GKW21b,Hay65,Hsu96} that the set $\Ecal^{\ell,t}$ of all equivalence classes forms an abelian group  under the multiplication
\[\langle f\rangle \langle g\rangle=\langle fg\rangle.
 \]
(When $t>0$, it is assumed that the constant term is nonzero.) It is also easy to see \cite{GKW21b,Hsu96} that
\[
\left|\Ecal^{\ell,t}\right|=(q-\llbracket t>0\rrbracket)q^{\ell+t-1}.
\]
For typographical convenience, we shall use  $\Ecal^{\ell}$ to denote $\Ecal^{\ell,0}$. For any given $q$, the following observation \cite[Lemma~1.1]{Hsu96} will be useful:
\begin{align}\label{eq:Egroup}
\Ecal^{\ell,t}\cong \Ecal^{\ell}\times \Ecal^{t-1} \times \fqx,\quad t\ge 1.
\end{align}
Thus we may focus on the group $\Ecal^{\ell}$.
Since $\Ecal^{\ell}$ is abelian, it is isomorphic to a direct product of cyclic groups.

 Let $\xi_{\ell,1},\ldots,\xi_{\ell,u_{\ell}}$ be a fixed set of generators of $\Ecal^{\ell}$, and denote their orders by $r_{\ell,1},\ldots,r_{\ell,u_{\ell}}$, respectively. In the rest of the paper, $\gamma$  denotes a fixed generator of $\fqx$. Each $\eps\in \Ecal^{\ell,t}$ can be written uniquely as
\[
\eps=\gamma^{e_0(\eps)}\prod_{h=1}^{u_{\ell}}\xi_{\ell,h}^{e_{\ell,h}(\eps)}
\prod_{i=1}^{u_{t-1}}\xi_{t-1,i}^{e_{\ell,i}(\eps)}.
\]
We shall call $\vec{e}(\eps)=(e_0(\eps),e_{\ell,1}(\eps),\ldots e_{\ell.u_{\ell}},e_{t-1,1},\ldots, e_{t-1,u_{t-1}})$ the {\em exponent vector} of $\eps$.  When $t=0$, it is understood that $e_0$ is ignored.
In the rest of the paper, we shall use $I_q(d;\eps)$ to denote the number of  polynomials in $\Ical_q(d)$ with leading and ending coefficients prescribed by $\eps\in\Ecal^{\ell,t}$, and $S_q(d;\eps)$ to denote the number of  polynomials in $\Scal_q(2d)$ with leading coefficients prescribed by $\eps\in \Ecal^{\ell}$.

Let $\om_r=\exp(2\pi i/r)$ and $\eps,\eps'\in \Ecal^{\ell,t}$.  Define
\begin{align}
\{\eps^{1/k}\}&=\{\delta\in  \Ecal^{\ell,t}:\delta^k=\eps\},\\
\Ecal^{\ell,t}(d)&=\{\langle f \rangle: f\in \Mcal_q(d)\}, \\
a(\eps,\eps')&=\om_{q-1}^{e_0(\eps)e_0(\eps')}\prod_{h=1}^{u_{\ell}}  \om_{r_{\ell,h}}^{e_{\ell,h}(\eps)e_{\ell,h}(\eps')}\prod_{i=1}^{u_{t-1}}  \om_{r_{t-1,i}}^{e_{t-1,i}(\eps)e_{t-1,i}(\eps')},\label{eq:a}\\
c(d;\eps) &= \sum_{\eps'\in \Ecal^{\ell,t}(d)}a(\eps,\eps'),\label{eq:cdj}\\
P(z;\eps)&=1+\sum_{d=1}^{\ell+t-1}c(d;\eps)z^d,\label{eq:Pj}\\
D&=\sum_{\eps\ne \langle 1\rangle}\deg(P(z;\eps)). \label{eq:D}
\end{align}
When $t>1$ and $\eps\in \Ecal^{\ell,t}$, we also write $\eps=(\eps_1,\gamma^n,\eps_2)$ with $\eps_1\in \Ecal^{\ell}$, $\eps_2\in \Ecal^{t-1}$, and $0\le n\le q-2$.
Define
\begin{align}\label{eq:D1}
D'&=\sum_{\delta\in \Ecal^{\ell}\setminus \{\langle 1\rangle\}}\deg(P(z;\delta,1,\delta)).
\end{align}
Since $\deg(P(z;\eps))\le \ell+t-1$ and $\deg(P(z;\delta,1,\delta))\le 2\ell$, we have

\begin{align}
D&\le (\ell+t-1)(\left|\Ecal^{\ell,t}\right|-1), \label{eq:Dbound1}\\
D'&\le 2\ell \left(q^{\ell}-1\right). \label{eq:Dbound2}
\end{align}
We shall use the Iverson bracket $\llbracket P  \rrbracket$ which has value 1 if the predicate $P$ is true and has value 0 otherwise.


Now we are ready to state our main results.
\begin{thm}\label{thm:thm1} Let $\Ecal$ denote the group $\Ecal^{\ell,t}$, $\eps\in \Ecal$ and  $D$ be defined in \eqref{eq:D}.  \\
(a) We have the following upper bound:\begin{align}
I_q(d;\eps)&\le \frac{1}{|\Ecal|}\frac{q^d-\llbracket t>0\rrbracket}{d} +\frac{D}{|\Ecal|}\frac{q^{d/2}}{d}.\label{eq:Iupper}
\end{align}
(b) Assume $\ell+t\le \lceil d/2\rceil -1$ and let $e_1(q,d)=\min\left\{3.4q^{-d/6}, 0.8\right\}$.  We have the following lower bound:
\begin{align}
I_q(d;\eps)&\ge\frac{1}{|\Ecal|}\frac{q^d-\llbracket t>0\rrbracket}{d}-
\left(\frac{D+\left|\{\eps^{1/2}\}\right|\llbracket 2\mid d\rrbracket}{|\Ecal|}+e_1(q,d)\right)\frac{q^{d/2}}{d}. \label{eq:IlowerAsy}
\end{align}
\end{thm}

 By \eqref{eq:cdj}--\eqref{eq:D}, we have $D\le (\ell+t-1)(|\Ecal|-1)$. Thus Theorem~1(a) immediately implies the following upper bound from \cite[Theorem~2.4]{Hsu96}.
\begin{cor}\label{cor:cor1}  We have
\begin{align}
I_q(d;\eps)&
\le \frac{1}{|\Ecal|}\frac{q^d}{d} +\frac{(|\Ecal|-1)(\ell+t-1)}{|\Ecal|}\frac{q^{d/2}}{d}.\label{eq:Hsu}
\end{align}
\end{cor}

Since $|\{\eps^{1/2}\}|\le |\Ecal|$, Theorem~1(b) immediately implies the following lower bound from \cite[Theorem~2.1]{Coh05}.
\begin{cor}\label{cor:cor2} Assume $\ell+t\le \lceil d/2\rceil -1$. Then
\begin{align}
I_q(d;\eps)&\ge \frac{1}{|\Ecal|}\frac{q^d}{d}-(\ell+t+1)\frac{q^{d/2}}{d}.\label{eq:Coh}
\end{align}
\end{cor}
We remark that the assumption $\ell+t\le \lceil d/2\rceil -1$ is without loss of generality because the right hand side of
 \eqref{eq:Coh} becomes negative when $\ell+t\ge  \lceil d/2\rceil$.  We also  note that the upper bound in \cite[Theorem~2.1]{Coh05}
is slightly weaker than \eqref{eq:Hsu}, and the lower bound in \cite[Theorem~2.4]{Hsu96} is slightly weaker than \eqref{eq:Coh}.

By \eqref{eq:Dbound1} and \eqref{eq:Dbound2}, our next theorem improves the error terms in \cite[Theorem~3]{Gao21} by a factor of $q^{\ell}$. This improvement
enables us to essentially extend the range of $\ell$ from $d/4$ to $d/2$.
\begin{thm}\label{thm:thm2} Let $\eps\in \Ecal^{\ell}$, $D$ be defined in \eqref{eq:D}, and $D'$ be defined in \eqref{eq:D1}. Assume $\ell\le \lceil d/2\rceil -1$ and let $e_2(q,d)=\min\left\{7q^{-d/6},2\right\}$. \\
(a) We have the following upper bound:
\begin{align}
 S_q(d;\eps)&\le \frac{1}{2d}q^{d-\ell}+\left(\frac{D'+2D+3\llbracket 2\mid d\rrbracket\left|\{\eps^{1/2}\}\right|}{2q^{\ell}}+e_2(q,d)\right)\frac{q^{d/2}}{d}.\label{eq:Supper}
\end{align}
(b)  We have the following lower bound:
\begin{align}
 S_q(d;\eps)&\ge \frac{1}{2d}q^{d-\ell}-\left(\frac{D'+2D}{2q^{\ell}}+\llbracket 2\mid d \rrbracket\frac{\left|\{\langle 1\rangle^{1/2}\}\right|}{q^{\ell}}+e_2(q,d)\right)\frac{q^{d/2}}{d}. \label{eq:Slower}
\end{align}
Consequently $S_q(d;\eps)>0$ whenever
\[
\ell\le \min\left\{\left\lceil \frac{d}{2} \right\rceil -1,\frac{d}{2}-\log_q (2d+2)\right\}.
\]
\end{thm}
\section{Proofs}

Set $g_0=1$. As in \cite{Gao21}, let $\phi_d:\Ecal^{\ell} \mapsto \Ecal^{\ell}$ be the bijection defined by
\begin{align}\label{eq:phi}
 f_k=\sum_{j\le k/2}{d+2j-k\choose j}g_{k-2j},~~1\le k\le \ell.
\end{align}

For each $\delta\in \Ecal^{\ell}$, we define the bijection
$\psi_{\delta}: \Ecal^{\ell}\mapsto \Ecal^{\ell}$ by
\begin{align}
\psi_{\delta}(\eps)=\eps\delta.\label{eq:psi}
\end{align}
We note that $\psi_{\delta}$ is the same as $\psi_{\vec{b}}$ defined in \cite{Gao21}, written in a more compact  notation.

The following  is \cite[Theorem~2]{Gao21}, rewritten in a slightly different notation.
\begin{thm} \label{thm:thm3} Suppose $d>1$ and $\eps\in \Ecal^{\ell}$. Then
\begin{align}
S_{q}\left(d;\eps\right)
&=\frac{1}{2}\sum_{\eps_1\in \{\eps^{1/2}\}} S_q(d/2;\eps_1) \nonumber\\
&~~~+I_q(d;\phi_d^{-1}(\eps))-\frac{1}{2}\sum_{n=0}^{q-2}\sum_{\delta\in \Ecal^{\ell}}
I_q\left(d;\eps\delta^{-1},\gamma^n,\delta\right).\label{eq:main}
\end{align}
\end{thm}

If $\{\eps^{1/k}\}\ne \emptyset$, we let $\eps^{1/k}$ to denote any particular element in $\{\eps^{1/k}\}$.
The following simple observations will be used later.
\begin{align}
\left\{\eps^{1/k}\right\}&=\eps^{1/k}\left\{\langle 1\rangle^{1/k}\right\}, \\
\left|\left\{\eps^{1/k}\right\}\right|&=\left|\left\{\langle 1\rangle^{1/k}\right\}\right|\left\llbracket \left\{\eps^{1/k}\right\}\ne \emptyset \right\rrbracket ,\\
a(\eps^{-1},\delta)&=a(\eps,\delta^{-1}),\label{eq:a1}\\
a(\delta,\eps_1\eps_2)&=a(\delta,\eps_1)a(\delta,\eps_2).\label{eq:a2}
\end{align}

Following Granger's notation \cite{Gra19},  we set
\[
\rho_d(g):=\sum_{\rho} \rho^{-d},
\]
where the sum is over all the nonzero roots (with  multiplicity) of the  polynomial $g\in {\mathbb C}[z]$.
Define
\begin{align}\label{eq:FI}
F_q(d;\eps)&=\sum_{k\mid d}\frac{d}{k}\sum_{\eps_1\in \{\eps^{1/k}\}} I_q(d/k;\eps_1).
\end{align}

Under the above notations, we can restate \cite[Theorem~3]{GKW21b} as follows.
\begin{thm}  \label{thm:thm4}
Let $\Ecal$ denote the group $\Ecal^{\ell,t}$ and let $\eps\in \Ecal$. We have
\begin{align}
I_q\left(d;\eps\right)&=\frac{1}{d}\sum_{k|d}\mu(k)\sum_{\eps_1\in \{\eps^{1/k}\}}
F_q(d/k;\eps_1),\label{eq:IF}\\
F_q(d;\eps)&=\frac{q^d-\llbracket t>0 \rrbracket}{|\Ecal|}+\frac{d}{|\Ecal|}\sum_{\delta\in \Ecal \setminus \{\langle 1\rangle\}}
a\left(\delta,\eps^{-1}\right)[z^d]\ln  P(z;\delta)\nonumber\\
 &=\frac{q^d-\llbracket t>0 \rrbracket}{|\Ecal|}-\frac{1}{|\Ecal|}\sum_{\delta\in \Ecal \setminus \{\langle 1\rangle\}}
a\left(\delta,\eps^{-1}\right)\rho_d(P(z;\delta)).\label{eq:Froot}
\end{align}
In particular, we have
\begin{align}
F_q(d;\langle 1\rangle)
 &= \frac{1}{|\Ecal|}\left(q^d-\llbracket t>0\rrbracket\right)-\frac{1}{|\Ecal|}
\sum_{\delta\in \Ecal \setminus \{\langle 1\rangle\}}\rho_d(P(z;\delta)). \label{eq:trivial}
\end{align}
\end{thm}

The following lemma simplifies sums involving $F_q(d;\eps)$ over some subgroups of $\Ecal^{\ell,t}$, which plays a crucial role in the proofs of Theorems~1 and 2.
\begin{lemma}\label{lemma1}  Let $\Ecal$ denote the group $\Ecal^{\ell,t}$ and $\Ecal^{\ell}$ denote the group $\Ecal^{\ell,0}$. \\
(a) For each $\eps\in \Ecal$, we have
\begin{align}
&~~~\sum_{\eps_1\in \{\eps^{1/k}\}}F_q\left(d/k;\eps_1\right)\nonumber\\
&= \frac{\left|\{\eps^{1/k}\}\right|}{|\Ecal|}\left(q^{d/k}-\llbracket t>0 \rrbracket\right)\label{eq:Fsum1}\\
&~~ -\frac{\left|\{\eps^{1/k}\}\right|}{|\Ecal|}\sum_{\delta\in \Ecal\setminus \{\langle 1\rangle\}} \left\llbracket \{\delta^{1/k}\}\ne \emptyset \right\rrbracket a\left(\delta,\eps^{-1/k}\right)\rho_{d/k}(P(z;\delta))\nonumber\\
&\le \frac{\left|\{\eps^{1/k}\}\right|}{|\Ecal|}\left(q^{d/k}-\llbracket t>0 \rrbracket\right)+(\ell+t-1) q^{d/2k}.\label{eq:Fsum2}
\end{align}
(b) Let $\gamma$ be any fixed generator of $\fqx$. For each $\eps\in \Ecal^{\ell}$, we have
\begin{align}
&~~~\sum_{m=0}^{q-2}\sum_{\delta\in \Ecal^{\ell}} F_q\left(d;\eps\delta^{-1},\gamma^m,\delta\right)\nonumber\\
&=\frac{q^d-1}{q^{\ell}}-\frac{1}{q^{\ell}}\sum_{\delta\ne \langle 1\rangle}a\left(\delta,\eps^{-1/k}\right)\rho_d\left(P(z;\delta,1,\delta)\right),
\label{eq:Fsum3}\\
&~~~\sum_{\eps_1\in \{\eps^{1/k}\}}\sum_{m=0}^{q-2}\sum_{\delta\in \Ecal^{\ell}} F_q\left(d/k;\eps_1\delta^{-1},\gamma^m,\delta\right)\nonumber\\
&=\frac{\left|\{\eps^{1/k}\}\right|}{q^{\ell}}\left(q^{d/k}-1-
\sum_{\delta\ne  \langle 1\rangle,\{\delta^{1/k}\}\ne \emptyset}a\left(\delta,\eps^{-1/k}\right)\rho_{d/k}\left(P(z;\delta,1,\delta)\right)\right).
\label{eq:Fsum4}
\end{align}
\end{lemma}
\proof The well-known identity
\[
\sum_{s=0}^{r-1}\om_{r}^{sj}=r\llbracket r\mid j\rrbracket
\]
immediately leads to
\begin{align}
\sum_{\delta\in \Ecal}a(\delta,\eps)&=|\Ecal|\llbracket \eps=\langle 1\rangle\rrbracket,\label{eq:aSum0}\\
\sum_{\eps_1\in \{\eps^{1/k}\}}a(\eps_1,\delta)&=\left|\left\{\eps^{1/k}\right\}\right|a(\eps^{1/k},\delta)\llbracket \{\delta^{1/k}\}\ne \emptyset\rrbracket. \label{eq:aSum}
\end{align}
where $a(\eps^{1/k},\delta)$ is interpreted as 0 if $\{\eps^{1/k}\}=\emptyset$.
It follows from  \eqref{eq:Froot} and \eqref{eq:aSum} that
\begin{align*}
&~~~\sum_{\eps_1\in \{\eps^{1/k}\}}F_q\left(d/k;\eps_1\right)\\
&=\frac{\left|\{\eps^{1/k}\}\right|}{|\Ecal|}\left(q^{d/k}-\llbracket t>0 \rrbracket\right)\\
&~~ -\frac{1}{|\Ecal|}\sum_{\delta\in \Ecal\setminus \{\langle 1\rangle\}}
\sum_{\eps_1\in \{\eps^{1/k}\}}a(\eps_1,\delta^{-1})\rho_{d/k}(P(z;\delta))\\
&=\frac{\left|\{\eps^{1/k}\}\right|}{|\Ecal|}\left(q^{d/k}-\llbracket t>0 \rrbracket\right)\\
&~~ -\frac{\left|\{\eps^{1/k}\}\right|}{|\Ecal|}\sum_{\delta\in \Ecal\setminus \{\langle 1\rangle\}} \left\llbracket \{\delta^{1/k}\}\ne \emptyset \right\rrbracket  a\left(\delta,\eps^{-1/k}\right)\rho_{d/k}(P(z;\delta)),
\end{align*}
which is \eqref{eq:Fsum1}. Now \eqref{eq:Fsum2} follows by noting
\begin{align} \label{eq:Esum}
\sum_{\delta\in \Ecal\setminus \{\langle 1\rangle\}} \left\llbracket \{\delta^{1/k}\}\ne \emptyset \right\rrbracket =\frac{|\Ecal|}{\left|\{\langle 1\rangle^{1/k}\}\right|}-1.
\end{align}
To prove part~(b), we use \eqref{eq:a1}, \eqref{eq:a2}, \eqref{eq:Froot}, and \eqref{eq:aSum0} to obtain
\begin{align*}
&~~~~\sum_{m=0}^{q-2}\sum_{\delta\in \Ecal^{\ell}} F_q\left(d;\eps\delta^{-1},\gamma^m,\delta\right)\nonumber\\
&=\frac{(q-1)q^{\ell}}{|\Ecal^{\ell,\ell+1}|}\left(q^d-1\right)\nonumber\\
&~~~-\frac{1}{|\Ecal^{\ell,\ell+1}|}\sum_{(\eps_1,\gamma^n, \eps_2)\ne \langle 1\rangle}a\left(\eps,\eps_1^{-1}\right)\rho_d\left(P(z;\eps_1,\gamma^n,\eps_2)\right) \sum_{m=0}^{q-2}\om_{q-1}^{-nm}\sum_{\delta\in \Ecal^{\ell}}a\left(\delta,\eps_1\eps_2^{-1}\right)\nonumber\\
&=\frac{1}{q^{\ell}}\left(q^d-1\right)-\frac{1}{q^{\ell}}\sum_{\delta\ne \langle 1\rangle} a(\eps^{-1},\delta)\rho_d\left(P(z;\delta,1,\delta)\right).
\end{align*}
Now \eqref{eq:Fsum4} follows from \eqref{eq:aSum}
\qed

\medskip

\noindent {\em Proof of Theorem~\ref{thm:thm1}}:
Hsu \cite[Theoem~1.3]{Hsu96} showed that each (complex) root $\rho$ of $P(z;\eps)$ satisfies the property
\begin{align} \label{eq:Weil}
\rho= 1~~\hbox{ or } ~~|\rho|=q^{-1/2}.
\end{align}
It follows from \eqref{eq:Froot} that
\begin{align}
\left|F_q(d;\eps)-\frac{1}{|\Ecal|}\left(q^d-\llbracket t>0\rrbracket \right)\right|
&\le \frac{D}{|\Ecal|}q^{d/2}.\label{eq:Fbound}
\end{align}
By \eqref{eq:FI}, we have
\begin{align}\label{eq:Flower}
F_q(d;\eps)&\ge dI_q(d;\eps)+\frac{d}{2}\sum_{\eps_1\in \{\eps^{1/2}\}}I_q(d/2;\eps_1).
\end{align}
It follows from \eqref{eq:Fbound} that
\begin{align}\label{eq:Iupper19}
I_q(d;\eps)&\le \frac{1}{d}F_q(d;\eps)
\le \frac{1}{d|\Ecal|}\left(q^d-\llbracket t>0\rrbracket\right) +\frac{D}{d|\Ecal|}q^{d/2},
\end{align}
which establishes the upper bound.

\medskip

We now prove the lower bound.
Define
\begin{align}
L_q(d,\ell,t)&=\sum_{k\ge 3}\llbracket k\mid d,\mu(k)=-1\rrbracket \frac{\left|\{\langle 1\rangle^{1/k}\}\right|}{|\Ecal|}q^{d/k-d/2}\label{eq:Lqd}\\
&~~~+(\ell+t-1)\sum_{k\ge 2}\llbracket k\mid d,\mu(k)=-1\rrbracket \left(1-\frac{\left|\{\langle 1\rangle^{1/k}\}\right|}{|\Ecal|} \right)q^{d/2k-d/2}.\nonumber
\end{align}

Using \eqref{eq:IF}, \eqref{eq:Fsum2}, and \eqref{eq:Fbound}, we obtain
\begin{align}
I_q(d;\eps)&\ge \frac{1}{d}F_q(d;\eps)-\frac{1}{d}\sum_{k\mid d}\llbracket \mu(k)=-1\rrbracket \sum_{\delta\in \{\eps^{1/k}\}}F_q(d/k;\delta) \nonumber \\
&\ge \frac{q^d-\llbracket t>0\rrbracket}{d|\Ecal|}
-\left(\frac{D+\left|\{\eps^{1/2}\}\right|\llbracket 2\mid d\rrbracket}{|\Ecal|}+L_q(d,\ell,t)\right)\frac{q^{d/2}}{d}.\label{eq:Ibound0}
\end{align}

 We now estimate $L_q(d,\ell,t)$ by truncating the sums in \eqref{eq:Lqd} and bounding the remainders by geometric sums. For our purpose, we use
\begin{align}
L_q(d;\ell,t)&\le  \sum_{k=3}^{29}\llbracket k\mid d,\mu(k)=-1\rrbracket \frac{\left|\{\langle 1\rangle^{1/k}\}\right|}{|\Ecal|}q^{d/k-d/2}\nonumber\\
&~~+(\ell+t-1)\left(1-\frac{1}{|\Ecal|}\right)\sum_{k=2}^{29}\llbracket k\mid d,\mu(k)=-1\rrbracket q^{d/2k-d/2}\label{eq:SumL}\\
&~~+\llbracket d\ge 30\rrbracket q^{-d/2}\sum_{1\le j\le d/30}\left(q^j+(\ell+t-1)\left(1-\frac{1}{|\Ecal|}\right)q^{j/2}\right).\nonumber
\end{align}
Using  \eqref{eq:SumL}, $\ell+t\le \lceil d/2\rceil -1$, and
 \[
 1-\frac{1}{|\Ecal|}<1,~\frac{1}{|\Ecal|}\le\frac{\left|\{\langle 1\rangle^{1/k}\}\right|}{|\Ecal|}\le 1,~\sum_{1\le j\le d/30}q^j\le \frac{q}{q-1} \left(q^{d/30}-1\right),
 \]
we obtain (using Maple) $L_q(d,\ell,t)\le \min\{2.8q^{-d/6},0.6\}$ when $q\ge 3$.

For the case $q=2$, we make use of the following observation: $|\Ecal|=2^{\ell+t-1}$ and
 $|\{\langle 1\rangle^{1/k}\}|=1$ when $k$ is an odd prime. Consequently we obtain from
\eqref{eq:SumL} that
\begin{align*}
L_2(d;\ell,t)&\le  2^{1-\ell-t}\sum_{k=3}^{29}\llbracket k\mid d,\mu(k)=-1\rrbracket  2^{d/k-d/2}
\nonumber\\
&~~+(\ell+t-1)\left(1-2^{1-\ell-t}\right)\sum_{k=2}^{29}\llbracket k\mid d,\mu(k)=-1\rrbracket q^{d/2k-d/2}\\
&~~+\llbracket d\ge 30\rrbracket 2^{1-d/2}\left(2^{d/30}-1\right)\nonumber\\
&~~+\llbracket d\ge 30\rrbracket 2^{-d/2}(\ell+t-1)\left(1-2^{1-\ell-t}\right)
\frac{\sqrt{2}}{\sqrt{2}-1}\left(2^{d/60}-1\right)\nonumber\\
&\le \min\left\{3.4q^{-d/6},0.8\right\}.
\end{align*}
 \qed

\noindent {\bf Remark} Some comments are in order about the lower bound in \eqref{eq:IlowerAsy} and its proof.
\begin{itemize}
\item  We may use \eqref{eq:Flower} and the lower bound \eqref{eq:IlowerAsy} to reduce the upper bound in  \eqref{eq:Iupper}.
\item We did not optimize the numerical values appeared in the estimations of $L_q(d,\ell,t)$. Improvement can be made by observing that $|\{\eps^{1/k}\}|/|\Ecal|$ is usually much smaller than 1. For example, when $2\nmid q$, we have $\left|\{\langle 1\rangle^{1/2}\}\right|=1+\llbracket t>0\rrbracket$.
\end{itemize}

\medskip

\noindent {\em Proof of Theorem~\ref{thm:thm2}}:
We first use Theorem~\ref{thm:thm4} and \eqref{eq:Esum} to simplify the following sum:
\begin{align}
&~~\sum_{m=0}^{q-2}\sum_{\delta\in \Ecal}I_q(d;\eps\delta^{-1},\gamma^m,\delta)\nonumber\\
&=\sum_{k\mid d}\frac{\mu(k)}{d}\sum_{m,\delta}\sum_{m_1,\delta_1,\eps_1}F_q(d/k;\eps_1\delta_1^{-1},\gamma^{m_1},\delta_1)
\llbracket \gamma^{km_1}=\gamma^m,\delta_1^k=\delta,\eps_1^k=\eps \rrbracket\nonumber\\
&=\sum_{k\mid d}\frac{\mu(k)}{d}\sum_{\eps_1\in \{\eps^{1/k}\}} \sum_{m_1=0}^{q-2}\sum_{\delta_1\in \Ecal}F_q(d/k;\eps_1\delta_1^{-1},\gamma^{m_1},\delta_1).\label{eq:Psi0}
\end{align}

Applying Theorem~\ref{thm:thm3}, \eqref{eq:Psi0}, and noting
\[
S_q(d/2;\delta)\le I_q(d/2;\delta)\le \frac{2}{d}F_q(d/2;\delta),
\]
 we obtain
\begin{align}
S_q(d;\eps)&\le \frac{1}{d}\sum_{\eps_1\in \{\eps^{1/2}\}} F_q(d/2;\eps_1)+I_q(d;\phi_d^{-1}(\eps))\nonumber\\
&~~~ -\frac{1}{2d}\sum_{m=0}^{q-2}\sum_{\delta\in \Ecal}F_q(d;\eps\delta^{-1},\gamma^{m},\delta) \label{eq:S9}\\
&~~+\sum_{k\ge 2}\frac{\llbracket k\mid d,\mu(k)=-1\rrbracket}{2d}\sum_{\eps_1\in \{\eps^{1/k}\}} \sum_{m=0}^{q-2}\sum_{\delta\in \Ecal}F_q(d/k;\eps_1\delta^{-1},\gamma^{m},\delta). \nonumber
\end{align}

It follows from  \eqref{eq:Iupper},  \eqref{eq:Fsum4}, and \eqref{eq:aSum}  that
\begin{align}
S_q(d;\eps)&\le \llbracket 2\mid d \rrbracket \frac{\left|\{\eps^{1/2}\}\right|}{q^{\ell}}\frac{q^{d/2}}{d}
+\llbracket 2\mid d \rrbracket (\ell-1)\left(1-\frac{1}{q^{\ell}}\right)  \frac{q^{d/4}}{d}\nonumber\\
&~~+\frac{1}{dq^{\ell}}\left(q^d
+Dq^{d/2}\right)-\frac{1}{2dq^{\ell}}\left(q^d-1\right)+\frac{D'}{2dq^{\ell}} q^{d/2}\nonumber \\
&~~+\llbracket 2\mid d \rrbracket \frac{\left|\{\eps^{1/2}\}\right| }{2q^{\ell}}\frac{q^{d/2}}{d}+ \llbracket 2\mid d \rrbracket \ell\left(1-\frac{1}{q^{\ell}}\right)\frac{q^{d/4}}{d}\nonumber \\
&~~~+\frac{1}{2dq^{\ell}}\sum_{k\ge 3}\llbracket k\mid d,\mu(k)=-1\rrbracket \left|\{\eps^{1/k}\}\right|q^{d/k} \nonumber\\
&~~+\frac{\ell}{d}\sum_{k\ge 3}\llbracket k\mid d,\mu(k)=-1\rrbracket \left(1-\frac{1}{q^{\ell}}\right)q^{d/2k}\nonumber\\
&\le \frac{1}{d}q^{d-\ell}+\left(
\frac{D'+2D+3\llbracket 2\mid d \rrbracket \left|\{\eps^{1/2}\}\right|}{2q^{\ell}}+U_q(d;\ell)\right)\frac{q^{d/2}}{d},\nonumber
\end{align}
where
\begin{align}
U_q(d;\ell)&=\frac{1}{2}q^{-\ell-d/2}+\llbracket 2\mid d \rrbracket (2\ell-1)\left(1-\frac{1}{q^{\ell}}\right) q^{-d/4}\nonumber\\
&~~~+\frac{1}{2q^{\ell}}\sum_{k\ge 3}\llbracket k\mid d,\mu(k)=-1\rrbracket \left|\{\eps^{1/k}\}\right|q^{d/k-d/2} \nonumber\\
&~~+\ell\left(1-\frac{1}{q^{\ell}}\right)\sum_{k\ge 3}\llbracket k\mid d,\mu(k)=-1\rrbracket q^{d/2k-d/2}.\label{eq:U}
\end{align}
Simple calculations as in the proof of Theorem~1 give  $U_q(d;\ell)\le \min\{6.6q^{-d/6},1.5\}$.

For the lower bound, we obtain from Theorem~\ref{thm:thm3} and \eqref{eq:Psi0} that
\begin{align}
S_q(d;\eps)&\ge I_q(d;\phi_d^{-1}(\eps))\nonumber\\
&~~-\sum_{k\mid d}\frac{\llbracket k\mid d\rrbracket}{2d}\sum_{\eps_1\in \{\eps^{1/k}\}} \sum_{m=0}^{q-2}\sum_{\delta\in \Ecal^{\ell}}F_q(d/k;\eps_1\delta^{-1},\gamma^{m},\delta) \nonumber\\
&\ge \frac{1}{d}q^{d-\ell}-\left(\frac{D+\left|\{\langle 1\rangle^{1/2}\}\right|\llbracket 2\mid d\rrbracket}{q^{\ell}}+L_q(d;\ell,0)\right)\frac{q^{d/2}}{d}\nonumber\\
&~~-\frac{1}{2d} \sum_{m=0}^{q-2}\sum_{\delta\in \Ecal}F_q(d;\eps\delta^{-1},\gamma^{m},\delta)\label{eq:Slow0}\\
&~~-\frac{1}{2d}\sum_{k\ge 3}\llbracket k\mid d,\mu(k)=1\rrbracket\sum_{\eps_1\in \{\eps^{1/k}\}} \sum_{m=0}^{q-2}\sum_{\delta\in \Ecal}F_q(d/k;\eps_1\delta^{-1},\gamma^{m},\delta).\nonumber
\end{align}

Applying \eqref{eq:Fsum3} and \eqref{eq:Esum} again, we obtain
\begin{align}
S_q(d;\eps)&\ge \frac{1}{d}q^{d-\ell}-\left(\frac{D+\left|\{\langle 1\rangle^{1/2}\}\right|\llbracket 2\mid d\rrbracket}{q^{\ell}}+L_q(d;\ell,0)\right)\frac{q^{d/2}}{d}\nonumber\\
&~~- \frac{1}{2dq^{\ell}}\left(q^d-1+D'q^{d/2}\right)\nonumber\\
&~~~-\frac{1}{2d}\sum_{k\ge 6}\llbracket k\mid d,\mu(k)=1 \rrbracket\frac{\left|\{\langle 1\rangle^{1/k}\}\right|}{q^{\ell}}q^{d/k}\nonumber\\
&~~~-\frac{\ell}{d}\left(1-\frac{1}{q^{\ell}}\right)\sum_{k\ge 6}\llbracket k\mid d,\mu(k)=1\rrbracket q^{d/2k}. \label{eq:Sbound}
\end{align}
Define
\begin{align}
L'_q(d;\ell)&=L_q(d;\ell,0)+\frac{1}{2}\sum_{k\ge 6}\llbracket k\mid d,\mu(k)=1 \rrbracket\frac{\left|\{\langle 1\rangle^{1/k}\}\right|}{q^{\ell}}q^{d/k-d/2}\nonumber\\
&~~~+\ell\left(1-\frac{1}{q^{\ell}}\right)\sum_{k\ge 6}\llbracket k\mid d,\mu(k)=1\rrbracket q^{d/2k-d/2}.
\end{align}
We then have
\begin{align}
S_q(d;\eps)&\ge \frac{1}{2d}q^{d-\ell}-\left(\frac{D'}{2q^{\ell}}+\frac{D+\left|\{\langle 1\rangle^{1/2}\}\right|\llbracket 2\mid d\rrbracket}{q^{\ell}}+L'_q(d;\ell)\right)\frac{q^{d/2}}{d}.
\end{align}
Similar calculations give
\[
L'_q(d;\ell)\le \min \{7q^{-d/6},2  \}.
\]
It follows that $S_q(d;\eps)>0$ when
\[
q^{d/2}>2(2\ell+2)q^{\ell}.
\]
Using $2\ell\le d-1$ and taking $\log_q$ on both sides, we complete the proof.
\qed

\section{Examples}

In this section, we use some examples to demonstrate that  $\frac{D}{|\Ecal^{\ell}|-1}$ is smaller than $\ell-1$. Let $P(z;\eps)$ be defined in \eqref{eq:Pj} and define
\begin{align*}
d_j=|\{\eps\in \Ecal^{\ell}\setminus \{\langle 1\rangle\}: \deg(P(z;\eps))=j\}|,~~\vec{d}=(d_1,d_2,\ldots d_{\ell-1}).
\end{align*}
We note
\[
D=\sum_{j=1}^{\ell-1} jd_j.
\]

\begin{example}  Consider $q=2$ and $\ell=4$. From \cite[Example~4]{GKW21b}, we have
\begin{align*}
\vec{d}&=(2,4,8), \\
D&=2+2\times 4+3\times 8=34,\\
\frac{D}{|\Ecal^4|-1}&=\frac{34}{2^4-1}<2.3.
\end{align*}
\end{example}

\begin{example}  Consider $q=2$ and $\ell=5$. From \cite[Example~6]{GKW21b}, we have
\begin{align*}
\vec{d}&=(2,4,8,16), \\
D&=2+2\times 4+3\times 8+4\times 16=98,\\
\frac{D}{|\Ecal^5|-1}&=\frac{98}{2^5-1}<3.2.
\end{align*}
\end{example}

\begin{example}  Consider $q=3$ and $\ell=3$. From \cite[Example~5]{GKW21b}, we have
\begin{align*}
\vec{d}&=(6,18), \\
D&=6+2\times 18=42,\\
\frac{D}{|\Ecal^3|-1}&=\frac{42}{3^3-1}<1.62.
\end{align*}
\end{example}

We use the following result from \cite[Lemma~1]{GKW21b} to produce a few more examples.
\begin{lemma}\label{lemma:generator}
Let $q=p$ be a prime number. The generators of $\Ecal^{\ell}$ are
\begin{align*}
\left\{\langle x^j +1 \rangle: \gcd(p, j) =1, 1\le j \le \ell\right\},
\end{align*}
and the order of  $\langle x^j +1 \rangle$ is  equal to $p^{s_j}$, where
 $s_j$ is the smallest  positive integer such that $jp^{s_j} > \ell$.
\end{lemma}

\begin{example}  Consider $q=2$ and $\ell=6$. By Lemma~\ref{lemma:generator}, the group $\Ecal^6$ is generated by $\langle x+1\rangle$, $\langle x^3+1\rangle$, $\langle x^5+1\rangle$, of orders 8,4,2, respectively. Using Maple, we find
\begin{align*}
\vec{d}&=(2,4,8,16,32), \\
D&=2+2\times 4+3\times 8+4\times 16+5\times 32=258,\\
\frac{D}{|\Ecal^6|-1}&=\frac{258}{2^6-1}<4.1.
\end{align*}
\end{example}

\begin{example}  Consider $q=2$ and $\ell=7$. By Lemma~\ref{lemma:generator}, the group $\Ecal^7$ is generated by $\langle x+1\rangle$, $\langle x^3+1\rangle$, $\langle x^5+1\rangle$, $\langle x^7+1\rangle$, of orders 8,4,2,2, respectively. Using Maple, we find
\begin{align*}
\vec{d}&=(2,4,8,16,32,64), \\
D&=2+2\times 4+3\times 8+4\times 16+5\times 32+6\times 64=642,\\
\frac{D}{|\Ecal^7|-1}&=\frac{642}{2^7-1}<5.1.
\end{align*}
\end{example}

\begin{example}  Consider $q=2$ and $\ell=8$. By Lemma~\ref{lemma:generator}, the group $\Ecal^8$ is generated by $\langle x+1\rangle$, $\langle x^3+1\rangle$, $\langle x^5+1\rangle$, $\langle x^7+1\rangle$, of orders 16,4,2,2, respectively. Using Maple, we find
\begin{align*}
\vec{d}&=(2,4,8,16,32,64,128), \\
D&=2+2\times 4+3\times 8+4\times 16+5\times 32+6\times 64+7\times 128=1538,\\
\frac{D}{|\Ecal^8|-1}&=\frac{1538}{2^8-1}<6.1.
\end{align*}
\end{example}

\begin{example}  Consider $q=3$ and $\ell=4$. By Lemma~\ref{lemma:generator}, the group $\Ecal^4$ is generated by $\langle x+1\rangle$, $\langle x^2+1\rangle$, $\langle x^4+1\rangle$,  of orders 9,3,3, respectively. Using Maple, we find
\begin{align*}
\vec{d}&=(6,18,54), \\
D&=6+2\times 18+3\times 54=204,\\
\frac{D}{|\Ecal^4|-1}&=\frac{204}{3^4-1}<2.6.
\end{align*}
\end{example}

\begin{example}  Consider $q=3$ and $\ell=5$. By Lemma~\ref{lemma:generator}, the group $\Ecal^5$ is generated by $\langle x+1\rangle$, $\langle x^2+1\rangle$, $\langle x^4+1\rangle$, $\langle x^5+1\rangle$,  of orders 9,3,3,3, respectively. Using Maple, we find
\begin{align*}
\vec{d}&=(6,18,54,162), \\
D&=6+2\times 18+3\times 54+4\times 162=852,\\
\frac{D}{|\Ecal^5|-1}&=\frac{852}{3^5-1}<3.6.
\end{align*}
\end{example}

\begin{example}  Consider $q=3$ and $\ell=6$. By Lemma~\ref{lemma:generator}, the group $\Ecal^6$ is generated by $\langle x+1\rangle$, $\langle x^2+1\rangle$, $\langle x^4+1\rangle$, $\langle x^5+1\rangle$,  of orders 9,9,3,3, respectively. Using Maple, we find
\begin{align*}
\vec{d}&=(6,18,54,162,486), \\
D&=6+2\times 18+3\times 54+4\times 162+5\times 486=3282,\\
\frac{D}{|\Ecal^6|-1}&=\frac{3282}{3^6-1}<4.51.
\end{align*}
\end{example}

\section{Conclusion} \label{conclusion}
We derived new error bounds for the number of irreducible monic polynomials with prescribed leading and ending coefficients. These bounds improve the bounds in \cite{Coh05,Hsu96}. The new bounds are then used to obtain bounds for the number $S_q(d;\eps)$ of self-reciprocal irreducible monic polynomials with prescribed leading coefficients. The new lower bound for $S_q(d;\eps)$ significantly improves that in \cite{Gao21} and it implies  $S_q(d;\eps)>0$ when
\[
\ell\le \min\left\{\left\lceil \frac{d}{2} \right\rceil -1,\frac{d}{2}-\log_q (2d+2)\right\},
\]
where $\eps$ can be any $\ell$ prescribed leading coefficients. Some examples are given to demonstrate the improvement of our bounds in Theorem~1 over those in \cite{Coh05,Hsu96}. Our examples show a pattern about the degree sequence $\vec{d}$, which may be used to calculate $D$. It would be interesting to see if $D/|\Ecal|$ is substantially smaller than $\ell+t-1$ for some $q$ and large $\ell+t$.

\medskip

\begin{center}
Acknowledgement
\end{center}

I would like to thank Prof.~Wan for helpful suggestions which improve the presentation of the paper. I also would like to thank Simon Kuttner for noticing an error in Example~6.


\begin{thebibliography}{999}








\bibitem{Coh05}
S.~D.~Cohen, Explicit theorems on generator polynomials,
{\em Finite Fields Appl.} {\bf 11}
(2005), 337--357.


\bibitem{Gao21} Z.C. Gao, Counting self-reciprocal irreducible monic polynomials with prescribed coefficients over a finite field, arXiv:2109.09006.

\bibitem{GarKap12} T. Garefalakis and G. Kapetanakis, On the Hansen-Mullen conjecture for self-reciprocal
irreducible polynomials, {\em Finite Fields Appl.} {\bf 69} (2012), 832--841.




\bibitem{GKW21b} Z.~Gao, S.~Kuttner, and Q.~Wang, Counting irreducible polynomials with prescribed coefficients over a finite field, arXiv:2109.02000.

\bibitem{Gra19} R. Granger, On the enumeration of irreducible polynomials over
$GF(q)$ with prescribed coefficients, {\em Finite Fields Appl.} {\bf 57} (2019), 156--229.

\bibitem{Ha16} J.~Ha, Irreducible polynomials with several prescribed coefficients, {\em Finite Fields Appl.} {\bf 40} (2016), 10--25.

\bibitem{Hay65} D. R. Hayes, The distribution of irreducibles in $GF[q, x]$, {\em Trans. Amer. Math. Soc.}
{\bf 117} (1965), 101--127.

\bibitem{Hsu96}
 C. N. Hsu, The distribution of irreducible polynomials in $\mathbb{F}_q[t]$,
  {\em J. Number Theory 61} {\bf (1)} (1996) 85--96.




\bibitem{PanTza12} D. Panario and G. Tzanakis, A generalization of the Hansen–Mullen conjecture on
irreducible polynomials over finite fields, {\em Finite Fields Appl.} {\bf 18} (2012) 303--315.


\bibitem{Pol13} P. Pollack, Irreducible polynomials with several prescribed coefficients, {\em Finite Fields Appl.} {\bf 22} (2013) 70--78.

\bibitem{Wan97} D. Wan, Generators and irreducible polynomials over finite fields, {\em Math. Comp.} {\bf 219} (1997) 1195--1212.
\end{thebibliography}
\end{document}